\documentclass[12pt]{article}
\usepackage{graphicx}
\usepackage{color}
\catcode`\@=11 \@addtoreset{equation}{section}

\catcode`\@=12

\newtheorem{Theorem}{Theorem}[section]
\newtheorem{Proposition}{Proposition}[section]
\newtheorem{Lemma}{Lemma}[section]
\newtheorem{Corollary}{Corollary}[section]
\newtheorem{Remark}{Remark}[section]
\newtheorem{Definition}{Definition}[section]

\newcommand{\bTheorem}[1]{
\begin{Theorem} \label{T#1} }
\newcommand{\eT}{\end{Theorem}}

\newcommand{\bProposition}[1]{
\begin{Proposition} \label{P#1}}
\newcommand{\eP}{\end{Proposition}}

\newcommand{\bLemma}[1]{
\begin{Lemma} \label{L#1} }
\newcommand{\eL}{\end{Lemma}}

\newcommand{\bCorollary}[1]{
\begin{Corollary} \label{C#1} }
\newcommand{\eC}{\end{Corollary}}

\newcommand{\bRemark}[1]{
\begin{Remark} \label{R#1} }
\newcommand{\eR}{\end{Remark}}

\newcommand{\bDefinition}[1]{
\begin{Definition} \label{D#1} }
\newcommand{\eD}{\end{Definition}}

\newcommand{\bFormula}[1]{
\begin{equation} \label{#1}}
\newcommand{\eF}{\end{equation}}
\newcommand{\bin}{\ \mbox{bounded in}\ }

\newcommand{\Ov}[1]{\overline{#1}}

\newcommand{\DC}{C^\infty_c}

\newcommand{\vr}{\varrho}
\newcommand{\vre}{\vr_\ep}

\newcommand{\vte}{\vt_\ep}
\newcommand{\vue}{\vu_\ep}
\newcommand{\tvr}{\tilde \vr}
\newcommand{\tvu}{{\tilde \vu}}
\newcommand{\tvt}{\tilde \vt}

\newcommand{\vt}{\vartheta}
\newcommand{\vu}{\vc{u}}
\newcommand{\vc}[1]{{\bf #1}}

\newcommand{\Div}{{\rm div}_x}
\newcommand{\Grad}{\nabla_x}

\newcommand{\tn}[1]{\mbox {\F #1}}
\newcommand{\dx}{{\rm d} {x}}
\newcommand{\dt}{{\rm d} t }

\newcommand{\intO}[1]{\int_{\Omega} #1 \ \dx}

\newcommand{\ep}{\varepsilon}

\font\F=msbm10 scaled 1000

\definecolor{Cgrey}{rgb}{0.85,0.85,0.85}
\definecolor{Cblue}{rgb}{0.50,0.85,0.85}
\definecolor{Cred}{rgb}{1,0,0}
\definecolor{fancy}{rgb}{0.10,0.85,0.10}

\newcommand\Cbox[2]{%
    \newbox\contentbox%
    \newbox\bkgdbox%
    \setbox\contentbox\hbox to \hsize{%
        \vtop{
            \kern\columnsep
            \hbox to \hsize{%
                \kern\columnsep%
                \advance\hsize by -2\columnsep%
                \setlength{\textwidth}{\hsize}%
                \vbox{
                    \parskip=\baselineskip
                    \parindent=0bp
                    #2
                }%
                \kern\columnsep%
            }%
            \kern\columnsep%
        }%
    }%
    \setbox\bkgdbox\vbox{
        \color{#1}
        \hrule width  \wd\contentbox %
               height \ht\contentbox %
               depth  \dp\contentbox
        \color{black}
    }%
    \wd\bkgdbox=0bp%
    \vbox{\hbox to \hsize{\box\bkgdbox\box\contentbox}}%
    \vskip\baselineskip%
}

\date{}

\begin{document}


\title{On the motion of viscous, compressible and heat-conducting liquids}

\author{Eduard Feireisl \thanks{The research of E.F. leading to these results has received funding from the European Research Council under the European Union's Seventh Framework
Programme (FP7/2007-2013)/ ERC Grant Agreement 320078. The Institute of Mathematics of the Academy of Sciences of
        the Czech Republic is supported by RVO:67985840.}\and Anton\' \i n Novotn\' y \and Yongzhong Sun}

\maketitle

\bigskip

\centerline{Institute of Mathematics of the Academy of Sciences of the Czech Republic}

\centerline{\v Zitn\' a 25, CZ-115 67 Praha 1, Czech Republic}

\bigskip

\centerline{Institut Math\'ematiques de Toulon, EA2134, University of Toulon}

\centerline{BP 20132, 839 57 La Garde, France }

\bigskip

\centerline{Department of Mathematics, Nanjing University}

\centerline{Nanjing, Jiangsu, 210093, P. R. China}






\maketitle

\bigskip





\begin{abstract}

We consider a system of equations governing the motion of a viscous, compressible, and heat conducting liquid-like fluid, with a general EOS of
Mie�-Gr{\" u}neisen type. In addition, we suppose that the viscosity coefficients may decay to zero for large values of the temperature.
We show the existence of global-in-time weak solution, derive a relative energy inequality, and compare the weak solutions with strong one
emanating from the same initial data - the weak strong uniqueness property.

\end{abstract}


\tableofcontents

\section{Introduction}
\label{i}

The motion of a viscous, compressible, and heat conducting fluid in continuum mechanics is described by a system of partial differential equations reflecting the three basic physical principles: The balance (or conservation) of mass, momentum, and energy expressed in terms of the mass density
$\vr = \vr(t,x)$, the macroscopic velocity $\vu(t,x)$, and the absolute temperature $\vt(t,x)$:
\bFormula{i1}
\partial_t \vr + \Div (\vr \vu) = 0,
\eF
\bFormula{i2}
\partial_t(\vr \vu) + \Div (\vr \vu \otimes \vu) + \Grad p(\vr, \vt) = \Div \tn{S}(\vt, \Grad \vu) + \vr \vc{f},
\eF
\bFormula{i3}
\partial_t \left( \frac{1}{2} \vr |\vu|^2 + \vr e(\vr, \vt) \right) +
\Div \left[ \left( \frac{1}{2} \vr |\vu|^2 + \vr e(\vr, \vt) + p(\vr, \vu) \right) \vu \right]
\eF
\[
+ \Div \vc{q}(\vt, \Grad \vt)
- \Div \left( \tn{S} (\vt, \Grad \vu ) \cdot \vu \right) = \vr \vc{f} \cdot \vu,
\]
where $p = p(\vr, \vt)$ is the pressure, $e = e(\vr, \vt)$ the internal energy, $\tn{S} = \tn{S}(\vt, \Grad \vu)$ the viscous stress tensor,
$\vc{q} = \vc{q}(\vt, \Grad \vt)$ the heat flux, and $\vc{f}$ the external driving force.

If the system is supplemented by conservative boundary conditions, say
\bFormula{i4}
\vu|_{\partial \Omega} = 0,\ \vc{q} \cdot \vc{n}|_{\partial \Omega} = 0,
\eF
equation (\ref{i3}) can be integrated over the physical domain $\Omega \subset R^3$ to obtain the total energy balance
\bFormula{i5}
\frac{{\rm d}}{{\rm d}t} \intO{ \left( \frac{1}{2} \vr |\vu|^2 + \vr e(\vr, \vt) \right) } =
\intO{\vr \vc{f} \cdot \vu}.
\eF

The thermodynamic functions $p = p(\vr, \vt)$ and $e = e(\vr, \vt)$ are interrelated through Gibbs' equation
\bFormula{i6}
\vt D s (\vr, \vt) = D e(\vr, \vt) + p(\vr, \vt) D \left( \frac{1}{\vr} \right),
\eF
where the new quantity $s = s(\vr, \vt)$ is called \emph{entropy}. Alternatively, the energy balance equation (\ref{i3}) can be written as
and internal energy equation
\bFormula{i7}
\partial_t \left( \vr e(\vr, \vt) \right) + \Div \left[ \vr e(\vr, \vt) \vu) \right] + \Div \vc{q}(\vt, \Grad \vt) = \tn{S}(\vt, \Grad \vu): \Grad \vu - p(\vr, \vt) \Div \vu
\eF
or the entropy equation
\bFormula{i8}
\partial_t \left( \vr s(\vr, \vt) \right) + \Div \left[ \vr s(\vr, \vt) \vu) \right] + \Div \left( \frac{ \vc{q} }{\vt} \right) =
\frac{1}{\vt} \left( \tn{S}(\vt, \Grad \vu):\Grad \vu - \frac{ \vc{q}(\vt, \Grad \vt) \cdot \Grad \vt}{\vt} \right),
\eF
with the entropy production rate
\[
\sigma = \frac{1}{\vt} \left( \tn{S}(\vt, \Grad \vu):\Grad \vu - \frac{ \vc{q}(\vt, \Grad \vt) \cdot \Grad \vt}{\vt} \right).
\]
As a consequence of the \emph{Second law of thermodynamics}, the entropy production rate $\sigma$ is non-negative for any physically admissible process - a stipulation that imposes certain restrictions on the constitutive relations for $\tn{S}$ and $\vc{q}$.

Our goal is to study the properties of system (\ref{i1}--\ref{i3}) provided the constitutive relations for $p$, $e$, $\tn{S}$, and $\vc{q}$ reflect the characteristic properties of \emph{liquids} rather than gases. In particular, we assume that:

\begin{itemize}

\item The pressure-density-temperature EOS is of Mie�-Gr{\" u}neisen type,
\bFormula{i9}
p_F(\vr, \vt) = p_c(\vr) + \vr \vt G(\vr, \vt),
\eF
where $p_c$ is the reference ``cold'' pressure, see \cite{EOS}, Shyue \cite{Shy} among others;
\item
The viscous stress tensor takes the Newtonian form
\bFormula{i10}
\tn{S}(\vt, \Grad \vu) = \mu(\vt) \left( \Grad \vu + \Grad^t \vu - \frac{2}{3} \Div \vu \tn{I} \right) + \eta (\vt) \Div \vu \tn{I},
\eF
where the transport coefficients $\mu$ and $\eta$ are bounded functions of the absolute temperature and may decay to zero for $\vt \to \infty$.

\end{itemize}

\bRemark{i1}

The fact that the viscosity coefficients may degenerate for $\vt \to \infty$ and the resulting
analytical difficulties are the main obstacles to be handled in the present paper.

\eR

Nonlinear systems like (\ref{i1}--\ref{i3}) are not known to possess global-in-time solution without imposing some extra restrictions on smallness and smoothness of the data. Alternatively, we can interpret the partial derivatives in the sense of distributions, and, following the seminal work of
Leray \cite{LER}, introduce the weak solutions. The theory of weak solutions for the full \emph{Navier-Stokes-Fourier} system
(\ref{i1}), (\ref{i2}), (\ref{i7}) or, alternatively, (\ref{i1}), (\ref{i2}), (\ref{i8}) was developed in \cite{EF70} and \cite{FENO6}, related results were obtained by Bresch and Desjardins \cite{BRDE}, \cite{BRDE1}, Hoff \cite{HOF6}, Lions \cite[Chapter 8, Section 8.5]{LI4}, among others.

All the afore-mentioned results concern \emph{gases} rather than \emph{liquids}, in particular, the viscosity coefficients are bounded below
away from zero, and even unbounded for $\vt \to \infty$ as in \cite{FENO6}. Here we focus on \emph{liquids}, the properties of which will be reflected in our choice of the EOS as well as the transport coefficients. Revisiting the abstract theory developed in \cite{FeNoSIMA}, we first show that the problem admits global-in-time weak solutions for any finite energy initial data, see Section \ref{g}. It is important that the associated weak formulation includes the entropy balance (\ref{i8}) as the latter gives rise to the associated \emph{relative energy inequality} that is crucial for the
existence of the so-called dissipative solutions. On the basis of this observation, we show the weak-strong uniqueness principle in Section \ref{ws}.
Finally, in Section \ref{cr}, we derive a conditional regularity result in the class of weak solutions.

\section{Preliminaries, hypotheses}
\label{P}

One of the principal and sofar unsurmountable difficulties in the theory of compressible fluid is the (hypothetical) possibility of formation of vacuum zones
- vanishing mass viscosity component of a weak solution - even if the initial density is strictly positive. As the quantities in the time derivatives in (\ref{i2}), (\ref{i7}) and/or (\ref{i8}) are always multiples of $\varrho$, the fields $\vu$, $e$ and/or $s$ may experience uncontrollable oscillations
in the vacuum. Similarly to \cite{FENO6}, we use the regularization effect of thermal radiation to control the behavior of $\vt$ on the vacuum zone. Accordingly, the pressure as well as the heat conduction are augmented by radiation components. Specifically, we consider the pressure
\bFormula{P1}
p(\vr, \vt) = p_F(\vr, \vt) + p_R(\vt), \ p_R(\vt) = \frac{a}{3} \vt^4, \ a > 0,
\eF
and the heat flux vector $\vc{q}$ in the form
\bFormula{P2}
\vc{q}(\vt, \Grad \vt) = - \kappa(\vt) \Grad \vt,
\eF
where
\bFormula{P3}
0 < \underline{\kappa} (1 + \vt)^\alpha \leq \kappa (\vt) \leq \Ov{\kappa} (1 + \vt)^\alpha,\ \alpha \geq 3,
\eF
cf. \cite[Chapter 1]{FENO6}.

\subsection{Constitutive relations}

In accordance with (\ref{P1}) and (\ref{i6}), we write
\bFormula{P4}
e(\vr, \vt) = e_F (\vr, \vt) + e_R (\vr, \vt), \ e_R = a \frac{\vt^4}{\vr},
\eF
where $p_F$, $e_F$ satisfy Gibbs' relation (\ref{i6}) for a certain entropy $s_F$. Moreover, we impose the hypothesis of thermodynamic stability:
\bFormula{P5}
\frac{\partial p_F(\vr, \vt)}{\partial \vr} > 0,\ \frac{\partial e_F(\vr, \vt)}{\partial \vt} > 0 \ \mbox{for all}\ \vr, \vt > 0.
\eF

Motivated by (\ref{i9}) and the abstract theory developed in \cite{FeNoSIMA}, we suppose
$p \in C^2((0, \infty)^2) \cap C^1([0, \infty)^2)$,
\bFormula{P6}
\lim_{\vr \to 0+} p_F(\vr, \vt) = 0 \ \mbox{for any}\ \vt > 0, \
\lim_{\vt \to 0+} p_F(\vr, \vt) = p_c(\vr) \ \mbox{for any}\ \vr > 0,
\eF
with the ``cold pressure'' $p_c$ satisfying
\bFormula{P7}
\underline{p} \vr^\gamma \leq p_c(\vr) \leq \Ov{p} (1 + \vr)^\gamma, \ \underline{p} > 0, \ \gamma > 3.
\eF
In addition, we suppose
\bFormula{P8}
\left| \frac{\partial p_F(\vr, \vt)}{\partial \vt} \right| \leq c \left( 1 + \vr^{\gamma/3} + \vt^3 \right)\ \mbox{for all}
\ 0 < \vt < \Theta_c (\vr),
\eF
where $\vr \mapsto \Theta_c(\vr)$ is a continuous curve satisfying
\bFormula{P9}
\Theta_c (\vr) \geq c \vr^{\gamma/4} - 1 \ \mbox{for a certain}\ c > 0.
\eF

As for the internal energy $e_F$, we assume
\bFormula{P10}
e_F (\vr, \vt) \geq 0,\  \lim_{[\vr, \vt] \to [0,0]} e_F(\vr, \vt) = 0,
\eF
\bFormula{P11}
c_v(\vr, \vt) \equiv \frac{\partial e_F(\vr, \vt)}{\partial \vt} \in C([0, \infty)^2),
\eF
\bFormula{P12}
0 < \underline{c} (1 + \vt)^\omega \leq c_v(\vr, \vt) \leq \Ov{c} (1 + \vt)^\omega \ \mbox{for all}\ \vr, \vt > 0.
\eF

\subsection{Viscosity coefficients}

The viscous stress $\tn{S}(\vt, \Grad \vu)$ is given by Newton's rheological law (\ref{i10}), where
$\mu, \ \eta \in W^{1,\infty}[0, \infty)$,
\bFormula{P13}
0 < \underline{\mu} (1 + \vt)^{-\beta} \leq \mu(\vt) \leq \Ov{\mu} (1 + \vt)^{-\beta} ,
\eF
\bFormula{P14}
0 < \underline{\eta} (1 + \vt)^{-\beta} \leq \eta(\vt) \leq \Ov{\eta} (1 + \vt)^{-\beta}
\eF
for all $\vt \in [0, \infty)$ and a certain $\beta \geq 0$.

\bRemark{P1}

In contrast with Newton's hypothesis, we suppose that the bulk viscosity coefficient $\eta$ is strictly positive. This assumption can be relaxed in the
case $\beta = 0$, for which the lower bound in (\ref{P14}) can be replaced by $\eta \geq 0$.

\eR

\section{Existence theory}
\label{g}

In this section, we introduced the concept of \emph{weak solution} to the Navier-Stokes-Fourier system and state out main existence result. To this end, we specify the initial conditions
\bFormula{g1}
\vr(0, \cdot) = \vr_0,\ \vt(0, \cdot) = \vt_0, \ \vr_0, \vt_0 > 0 \ \mbox{in} \ \Omega,\ \vu(0, \cdot) = \vu_0.
\eF

\subsection{Weak solutions}
\label{WS}

Following \cite[Chapter 2]{FENO6} we introduce the weak solutions to problem (\ref{i1}--\ref{i4}), (\ref{g1}).

\subsubsection{Equation of continuity - mass conservation}

Equation (\ref{i1}) is understood in the sense of weak renormalized solutions introduced by DiPerna and Lions \cite{DL}:
\bFormula{g2}
\intO{ \left[ \Big(\vr + b(\vr) \Big) \varphi (\tau, \cdot) - \Big( \vr_0 + b(\vr_0) \Big)  \varphi (0, \cdot) \right] }
\eF
\[
= \int_0^\tau \intO{ \left[ \Big( \vr + b(\vr) \Big) \partial_t \varphi + \Big( \vr + b(\vr) \Big) \vu \cdot \Grad \varphi +
\Big( b(\vr) - b'(\vr) \vr \Big) \Div \vu \varphi \right] } \ \dt
\]
for any $\tau \in [0,T]$, and $b \in C^1[0,\infty)$, $b'(\vr) = 0$ for $\vr>>1$, and any test function
$\varphi \in \DC([0,T] \times \Ov{\Omega})$.

\subsubsection{Momentum balance}

Equation (\ref{i2}) is replaced by a family of integral identities
\bFormula{g3}
\intO{ \left[ \vr \vu \cdot \varphi (\tau, \cdot) - \vr_0 \vu_0  \cdot \varphi (0, \cdot) \right] }
\eF
\[
= \int_0^\tau \intO{ \left[ \vr \vu \cdot \partial_t \varphi + \vr \vu \otimes \vu : \Grad \varphi + p(\vr, \vt) \Div \varphi
- \tn{S}(\vt, \Grad \vu) : \Grad \varphi + \vr \vc{f} \cdot \varphi \right] } \ \dt
\]
for any $\tau \in [0,T]$ and any $\varphi \in \DC([0,T] \times \Omega; R^3)$.

\subsubsection{Energy balance}

A proper formulation of the energy balance is the cornerstone of the existence theory. Here, similarly to \cite[Chapter 2]{FENO6}, we postulate the entropy inequality
\bFormula{g4}
\intO{ \left[ \vr s(\vr, \vt) \varphi (\tau, \cdot) - \vr_0 s(\vr_0, \vt_0) \varphi (0, \cdot) \right] }
\eF
\[
\geq \int_0^\tau \intO{ \left[ \vr s(\vr, \vt) \partial_t \varphi + \vr s(\vr, \vt) \vu \cdot \Grad \varphi + \frac{\vc{q}(\vt, \Grad \vt) }{\vt}
\cdot \Grad \varphi \right] } \ \dt
\]
\[
+ \int_0^\tau \intO{ \frac{1}{\vt} \left[ \tn{S}(\vt, \Grad \vu) + \frac{ \vc{q}(\vt, \Grad \vt) \cdot \Grad \vt }{\vt} \right] \varphi } \ \dt
\]
for a.a. $\tau \in [0,T]$ and any test function $\varphi \in \DC([0,T] \times \Ov{\Omega})$, $\varphi \geq 0$; together with the total energy balance
\bFormula{g5}
\intO{ \left[ \left( \frac{1}{2} \vr |\vu|^2 + \vr e(\vr, \vt) \right) \psi(\tau) - \left( \frac{1}{2} \vr_0 |\vu_0|^2 + \vr_0 e(\vr_0, \vt_0) \right) \psi(0) \right]  }
\eF
\[
= \int_0^\tau \left[ \intO{ \left( \frac{1}{2} \vr |\vu|^2 + \vr e(\vr, \vt) \right) } \right] \partial_t \psi \ \dt
+ \int_0^\tau \intO{ \vr \vc{f} \cdot \vu } \psi \ \dt
\]
for a.a. $\tau \in [0,T]$ and any $\psi \in \DC[0, T]$.

In addition to \cite[Chapter 2]{FENO6}, as a benefit of the method used in the present paper, we are able to establish
the thermal energy balance
\bFormula{g6}
\intO{ \left[  \vr \Big( e(\vr, \vt) - e_c(\vr) \Big)   \psi(\tau) -  \vr_0 \Big( e(\vr_0, \vt_0)
- e_c (\vr_0) \Big)  \psi(0) \right]  }
\eF
\[
\geq \int_0^\tau \intO{ \vr (e - e_c) } \partial_t \psi \ \dt + \int_0^\tau \intO{ \left[ \tn{S}(\vt, \Grad \vu) : \Grad \vu -
(p - p_c) \Div \vu \right] } \psi \ \dt
\]
for a.a. $\tau \in [0,T]$ and any $\psi \in \DC[0,T]$, where we have set
\[
e_c(\vr) = \vr \int_1^\vr \frac{p_c(z)}{z^2} \ {\rm d}z,\ p_c(\vr) = p(\vr, 0),
\]
cf. \cite[Section 2.2]{FeNoSIMA}.

\subsection{Global in time weak solutions}

We are ready to state our main result concerning \emph{existence} of global-in-time weak solutions.

\Cbox{Cgrey}{

\bTheorem{g1}
Let $\Omega \subset R^3$ be a bounded domain of class $C^{2 + \nu}$. Suppose that the thermodynamic functions $p$, $e$, satisfy the structural hypotheses (\ref{P1}), (\ref{P4}--\ref{P12}), where
\bFormula{g7}
\gamma > 3, \ 0 \leq \omega \leq \frac{1}{2},
\eF
and that $s$ is determined (modulo an additive constant) by Gibbs' relation (\ref{i6}).
Let $\tn{S}$ and $\vc{q}$ be given by (\ref{i10}), (\ref{P2}), with the transport coefficients satisfying (\ref{P3}), (\ref{P13}), (\ref{P14}),
where
\bFormula{g8}
0\le \beta\le 4,\;\;\alpha\ge \frac {16} 3+\beta.
\eF
Let
\[
\vc{f} \in L^\infty((0,T) \times \Omega; R^3)
\]
be a given external force.
Suppose the initial data (\ref{g1}) satisfying
\bFormula{g9}
\vr_0 , \ \vt_0 \in L^\infty(\Omega),\ \vu_0 \in L^2(\Omega; R^3),\ \vr_0 \geq \underline{\vr} > 0,\ \vt_0 \geq
\underline{\vt} > 0 \ \mbox{a.a. in}\ \Omega.
\eF

Then the Navier-Stokes-Fourier system admits a weak solution $[\vr, \vt, \vu]$ in the space-time cylinder $(0,T) \times \Omega$ in the sense
specified in Section \ref{WS}. The weak solution belongs to the class
\[
\vr \in C([0,T]; L^1(\Omega)) \cap C_{\rm weak}([0,T]; L^\gamma (\Omega)), \ \vr \geq 0 \ \mbox{a.a. in}\ (0,T) \times \Omega;
\]
\[
\vt \in L^\infty(0,T; L^4(\Omega)) \cap L^2(0,T; W^{1,2}(\Omega)) \cap L^\alpha (0,T; L^{3 \alpha}(\Omega)),\ \vt > 0 \ \mbox{a.a. in}\
\Omega;
\]
\[
\vu \in L^q(0,T; W^{1,q}_0(\Omega; R^3)) \ \mbox{for a certain}\ q > 1.
\]

\eT

}

The rest of this section is devoted to the proof of Theorem \ref{Tg1}. With the machinery developed in \cite[Chapter 3]{FENO6}, the proof reduces basically to verification suitable uniform bounds ({\it a priori} estimates) imposed on the solutions set by the constraints specified through the list of hypotheses in Theorem \ref{Tg1}. We point out that these estimates apply to \emph{any} weak solution of the problem enjoying the properties specified in Section \ref{WS}.

\subsubsection{Uniform bounds}
\label{UB}

{\bf Step 1: Energy estimates}

\medskip

In accordance with (\ref{g2}), the total mass $\intO{ \vr }$ is a constant of motion; whence
\bFormula{g10}
\vr \bin L^\infty(0,T; L^1(\Omega)).
\eF
Next, it follows from (\ref{g5}) and the hypotheses imposed on the data that the total energy
\[
E = \intO{ \left( \frac{1}{2} \vr |\vu|^2 + \vr e(\vr, \vt) \right) }
\]
remains bounded on any compact time interval $[0,T]$. In accordance with (\ref{P4}), (\ref{P10}) we immediately deduce
\bFormula{g11}
\vt \bin L^\infty(0,T; L^4(\Omega)),\ \vr |\vu|^2 \bin L^\infty(0,T; L^1(\Omega)).
\eF
Since
$$
0\le e_F(\vr,\vt)-e_F(\vr,0)=\Big(e_F(\vr,\vt)-e_F(1,\vt)-e_F(\vr,0)\Big) +e_F(1,\vt)
$$
we may write
\bFormula{g12}
e_F(\vr, \vt) = e_c (\vr) + \int_0^\vt c_v(1,s) \ {\rm d}s + \int_1^\vr \frac{1}{z^2} \left( p_F(z, \vt) - p_F(z,0) -
\vt \frac{\partial p_F (z, \vt) }{\partial \vt} \right) \ {\rm d}z
\eF
and deduce exactly as in \cite[Section 4.2]{FeNoSIMA},
\bFormula{g13}
\vr \vt^{1 + \omega} \bin L^\infty(0,T; L^1(\Omega)),
\eF
\bFormula{g14}
\vr \bin L^\infty(0,T; L^\gamma(\Omega)),
\eF
and
\bFormula{g15}
p(\vr, \vt) \bin L^\infty(0,T; L^1(\Omega)).
\eF

\medskip

\noindent
{\bf Step 2: Dissipation estimates - entropy}

\medskip

Taking $\varphi \equiv 1$ in the entropy inequality (\ref{g4}) we obtain
\bFormula{g17-}
\intO{ \vr s(\vr, \vt) (\tau, \cdot) } \geq  \intO{ \vr_0 s(\vr_0, \vt_0)  }
+ \int_0^\tau \intO{ \frac{1}{\vt} \left[ \tn{S}(\vt, \Grad \vu) + \frac{ \vc{q}(\vt, \Grad \vt) \cdot \Grad \vt }{\vt} \right]  } \ \dt,
\eF
where
\bFormula{g17}
s(\vr, \vt) = s_F(\vr, \vt) + s_R(\vr, \vt),\ s_R(\vr, \vt) = \frac{4}{3} a \frac{\vt^3}{\vr},
\eF
and
\[
s_F(\varrho, \vartheta) =
s_F(\varrho,1) + \int_1^{\vartheta} \frac{ c_v (\varrho, s)}{s} \
{\rm d}s,\ s_F(\varrho,1) = -
\int_1^{\varrho} \frac{\partial p_F (z, 1)}{
\partial \vartheta } \frac{1}{z^2} \ {{\rm d}z}.
\]

Repeating the arguments of \cite[Section 4.3]{FeNoSIMA} we have
\bFormula{g18}
| s_F(\varrho,1) | \leq c \Big( 1 + \varrho^{-1}  +
\varrho^{{\gamma \over 3} - 1} \Big) ,
\eF
together with
\bFormula{g19}
{c\vr[\ln \vt]^+}\le\vr
\Big[ \int_1^{\vartheta} \frac{ c_v (\varrho, s)}{s} \ {\rm d}s
\Big]^+
\eF
\[
\leq \left\{
\begin{array}{c}
c{\vr}(1 + \vartheta^{\omega}),\ (\mbox{if} \ \omega>0)\\
c{\vr}(1+[\log(\vartheta)]^+),\ (\mbox{if}\ \omega=0)
\end{array}
\right\} \ \mbox{for all}\ { \vr>0},\vartheta > 0.
\]

In view of the energy bounds already established in (\ref{g11}--\ref{g15}), the left-hand side of (\ref{g17-}) is bounded from above
in $[0,T]$ and we may infer
\bFormula{g20}
\vr s(\vr, \vt) \bin L^\infty(0,T; L^1(\Omega)),
\eF
\bFormula{g21}
\frac{\tn{S}(\vt, \Grad \vu) : \Grad \vc{u}}{\vartheta} ,  \
\frac{\kappa ( \vartheta)}{\vartheta^2} |\Grad
\vartheta |^2
\bin L^1((0,T) \times \Omega),
\eF
and
\bFormula{g22}
{ c\vr[\ln\vt]^-\le}
\varrho \Big[ \int_1^{\vartheta}  \frac{ c_v ( \varrho, s )}{s}
\ {\rm d}s \Big]^- { \le c\Big(\vr|s(\vr,\vt)| +\vr e(\vr,\vt)\Big)}.
\eF

Thus, making use of hypotheses (\ref{P3}), we may follow step by step the arguments of \cite[Section 4.3]{FeNoSIMA} to obtain
\bFormula{g23}
\vr \log (\vt) \bin L^\infty(0,T; L^1(\Omega)),\ \log(\vt) \bin L^2(0,T; W^{1,2}(\Omega)),
\eF
\bFormula{g24}
\vt^{\alpha/2} \bin L^2(0,T; W^{1,2}(\Omega),
\eF
where the last relation combined with the standard embedding $W^{1,2}(\Omega) \hookrightarrow L^6(\Omega)$ give rise to
\bFormula{g25}
\vt \bin L^\alpha(0,T; L^{3\alpha}(\Omega)).
\eF

\medskip

\noindent
{\bf Step 3: Dissipation estimates - internal energy}

\medskip

The thermal energy balance (\ref{g6}), together with the bounds (\ref{g13}), (\ref{g14}), yields
\bFormula{g26}
\int_0^\tau \intO{ \tn{S}(\vt, \Grad \vu) : \Grad \vu } \ \dt \leq c \left(1 + \int_0^\tau \intO{ \left| p(\vr, \vt) -
p_c(\vr) \right||\Div \vu| } \ \dt \right).
\eF

Arguing as in \cite[Section 4.3]{FeNoSIMA} we get
\bFormula{g27}
| p(\varrho, \vartheta) - p_c(\varrho) | \leq c ( 1 + \vartheta^4 +
\vartheta \varrho^{\gamma \over 3} )
\ \mbox{for all} \ \varrho > 0 ,\ \vartheta > 0.
\eF
Writing
\[
\left| p(\vr, \vt) -
p_c(\vr) \right||\Div \vu| = \left| p(\vr, \vt) -
p_c(\vr) \right| \vt^{\beta/2} {\vt^{-\beta/2}}|\Div \vu|
\]
we observe easily that (\ref{g26}) combined with hypotheses (\ref{P13}), (\ref{P14}) will yield the desired conclusion
\bFormula{g28}
\tn{S}(\vt, \Grad \vu) : \Grad \vu \bin L^1((0,T) \times \Omega)
\eF
as soon as we establish a bound for
\bFormula{g29}
( 1 + \vartheta^4 +
\vartheta \varrho^{\gamma \over 3} ) \vt^{\beta/2} \ \mbox{in}\ L^2((0,T) \times \Omega).
\eF
To see (\ref{g29}), we interpolate between (\ref{g11}), (\ref{g25}) to obtain
\bFormula{g30}
\vt \bin L^q(0,T; L^p(\Omega)) \ \mbox{whenever} \ p = \frac{12 \alpha}{4 \lambda - 3 \lambda \alpha + 3 \alpha},\ q \leq \frac{\alpha}{\lambda},\
\lambda \in [0,1].
\eF
Consequently, for $q = 8 + \beta$, we have $\lambda = \frac{\alpha}{8 + \beta}$ with the corresponding
\[
p = (8 + \beta) \frac{ 12 }{4 - 3\alpha + 3 (8 + \beta)} \geq 8 + \beta
\]
whenever $\alpha$, $\beta$ satisfy hypothesis (\ref{g8}). Thus
\[
\vt^{4 + \beta/2} \bin L^2((0,T) \times \Omega).
\]
The last term in (\ref{g29}) can be handled in a similar fashion seeing that, in view of { (\ref{g25})},
\[
\vr^{\gamma/3} \bin L^\infty(0,T; L^{3}(\Omega)),\ \vt^{1 + \beta/2} \bin L^2(0,T; L^6(\Omega)).
\]
We have shown (\ref{g28}), in particular,
\bFormula{g31}
(1 + \vt)^{-\beta/2} |\Grad \vu + \Grad \vu^t | = g \bin L^2((0,T) \times \Omega).
\eF

Finally, using once more the estimates (\ref{g11}), (\ref{g25}), we deduce from (\ref{g31}) that
\[
|\Grad \vu + \Grad \vu^t | \bin L^2(0,T; L^q(\Omega)),\ q = \frac{8}{\beta + 4} \ \mbox{if}\ \beta \leq 4,
\]
and
\[
\| \vt^{\beta/2} g \|_{L^{6\alpha/ (3\alpha + \beta)} (\Omega)} \leq \| g \|_{L^2(\Omega)} \| \vt^{\beta/2} \|_{L^{6\alpha/ \beta}(\Omega)} =
\| g \|_{L^2(\Omega)} \| \vt \|_{L^{3\alpha}(\Omega)}^{\beta/2};
\]
whence
\[
|\Grad \vu + \Grad \vu^t | \bin L^p(0,T; L^r(\Omega)),\ r = \frac{6 \alpha}{3 \alpha + \beta}, \ p = \frac{2 \alpha}{\alpha + \beta}.
\]
Applying the standard Korn's and Poincar\' e's inequalities, we conclude
\bFormula{g32}
\left\{
\begin{array}{c}
\vu \bin L^2(0,T; W^{1,q}_0(\Omega; R^3)),\ q = \frac{8}{\beta + 4} \ \mbox{if}\ \beta \leq 4, \\ \\
\vu \bin L^p(0,T); W^{1,r}_0 (\Omega; R^3)),\ r = \frac{6 \alpha}{3 \alpha + \beta}, \ p = \frac{2 \alpha}{\alpha + \beta}.
\end{array}
\right\}
\eF

\medskip

\noindent

{\bf Step 4: Pressure estimates}

\medskip

The final estimates concern the pressure $p$ sofar bounded only by (\ref{g15}). In order to improve integrability with respect to time, the standard
procedure based on the application of the so-called Bogovskii's operator is needed, see \cite[Section 4.5]{FeNoSIMA}, \cite[Chapter 2, Section 2.2.5]{FENO6}. To perform this step and also for construction of the weak solutions we need the convective terms to be belong to the Lebesgue space
$L^p$ for some $p > 1$. It turns out, see \cite[Section 5]{FeNoSIMA}, that the velocity field must satisfy, at least
\bFormula{g33}
\vu \bin L^q(0,T; L^s(\Omega; R^3)) \ \mbox{for certain}\ q > 1,\ s > 4.
\eF
Revisiting (\ref{g32}) we need $W^{1,r}(\Omega) \hookrightarrow L^4(\Omega)$, meaning
$$
\frac{18\alpha}{3(\alpha+\beta)}>4\;\mbox{i.e.}\; \alpha> 2\beta,
$$
which is in agreement with hypothesis (\ref{g8}). We conclude, exactly as in \cite[Section 4.5]{FeNoSIMA},
\bFormula{g34}
p(\vr, \vt)^r \bin L^1((0,T) \times \Omega) \ \mbox{for some}\ r > 1.
\eF

\subsubsection{Proof of Theorem \ref{Tg1}}

We have collected all the ingredients necessary for the proof of Theorem \ref{Tg1}. This can be carried over as follows.

\begin{itemize}

\item We consider the problem with the shear viscosity coefficient augmented, specifically,
\[
\mu_\ep (\vt) = \mu (\vt) + \ep (1 + \vt),\ \ep > 0.
\]
Adapting the construction performed in \cite[Chapter 3]{FENO6} we obtain a family of approximate solutions $\{ \vre, \vue, \vte \}_{\ep > 0}$.

\item It is easy to check that $\{ \vre, \vue, \vte \}_{\ep > 0}$ admits the uniform bounds established in Section \ref{UB}

\item
Having established all the necessary estimates,
we apply the compactness arguments of \cite[Section 5]{FeNoSIMA} to pass to the limit for $\ep \to 0$ in $\{ \vre, \vue, \vte \}_{\ep > 0}$.

\end{itemize}

We have proved Theorem \ref{Tg1}.

\section{Relative energy and weak strong uniqueness}
\label{ws}

We study stability properties of the class of weak solutions, the existence of which was proved in Theorem \ref{Tg1}. In particular, we recall the relative energy inequality and show the weak-strong uniqueness principle. For the sake of simplicity, we take
\[
\vc{f} \equiv 0.
\]

\subsection{Relative energy}

The relative entropy/energy and the associated concept of dissipative solutions was introduced in the pioneering paper by Dafermos \cite{Daf4}
and later developed by Leger and Vasseur \cite{LegVas}, Masmoudi \cite{MAS5}, Saint-Raymond \cite{SaiRay},  among others. In particular, the problem of weak-strong uniqueness for the compressible Navier-Stokes and the Navier-Stokes-Fourier system
were addressed by Germain \cite{Ger} and finally solved in \cite{FeiNov10}, \cite{FeNoJi}, \cite{FeNoSun}.

We recall the relative energy associated to the Navier-Stokes-Fourier system
\bFormula{ws1}
\mathcal{E} \left( \vr, \vt, \vu \Big| r, \Theta, \vc{U} \right) =
\intO{ \left[ \frac{1}{2} \vr |\vu - \vc{U} |^2 + H_\Theta (\vr, \vt) - \frac{\partial H_\Theta(r, \Theta) }{\partial \vr} (\vr - r) -
H_\Theta (r, \Theta) \right] },
\eF
where
\[
H_\Theta(\vr, \vt) = \vr \Big( e(\vr, \vt) - \Theta s(\vr, \vt) \Big).
\]
As shown in \cite[Chapter 3, Proposition 3.2]{FENO6}, the functional $\mathcal{E}$ represents a ``distance'' between the quantities
$[\vr, \vt, \vu]$ and $[r, \Theta, \vc{U}]$. Specifically, for any compact sets { $ K, K'$, $K'\subset {\rm int}(K) \subset (0, \infty)^2$}, there exists a positive constant
{ $c=c(K,K')$}, depending solely on the structural properties of the thermodynamic functions $e$ and $s$, notably on the lower bounds for the
derivatives appearing in the thermodynamics stability hypothesis (\ref{P5}), such that
\bFormula{ws2}
\mathcal{E} \left( \vr, \vt, \vu \Big| r, \Theta, \vc{U} \right) \geq c \left\{ \begin{array}{l}  |\vr - r|^2 + |\vt - \Theta|^2 + |\vu - \vc{U}|^2
\ \mbox{if}\ [\vr, \vt] \in K, \ [r, \Theta] \in { K'} \\ \\
1 + \vr |\vu - \vc{U}|^2 + \vr e(\vr, \vt) + \vr |s(\vr, \vt)|  \ \mbox{if} \ [\vr, \vt] \in (0, \infty)^2 \setminus K, \ [r, \Theta] \in { K'}.
\end{array} \right.
\eF

Any weak solution $[\vr, \vu, \vt]$ of the Navier-Stokes-Fourier system in the sense of Section \ref{WS} that belongs to the regularity class
specified in Theorem \ref{Tg1} satisfies the relative energy inequality
\bFormula{ws3}
\left[ \mathcal{E} \left( \vr , \vt , \vu \Big| r, \Theta, \vc{U} \right) \right]_{t = 0}^{t = \tau}
\eF
\[
 +
\int_0^\tau \intO{ \frac{\Theta}{\vt} \left(
\tn{S}(\vt, \Grad \vu): \Grad \vu - \frac{ \vc{q}(\vt, \Grad \vt) \cdot \Grad \vt }{\vt} \right) } \ \dt
\]
\[
\leq
\int_0^\tau \intO{ \vr (\vu - \vc{U}) \cdot \Grad \vc{U} \cdot (\vc{U} - \vu)} \ \dt
\]
\[
+ \int_0^\tau \intO{  \tn{S}(\vt, \Grad \vu) : \Grad \vc{U} } \ \dt - \int_0^\tau \intO{
\frac{ \vc{q}(\vt, \Grad \vt) }{\vt} \cdot \Grad \Theta  } \ \dt +  \lambda \int_0^\tau \intO{ \vu \cdot \vc{U} } \ \dt
\]
\[
+ \int_0^\tau \intO{ \vr \Big( s(\vr, \vt) - s(r, \Theta) \Big) \Big( \vc{U} - \vu \Big)
\cdot \Grad \Theta } \ \dt
\]
\[
+ \int_0^\tau \intO{ \vr \Big(  \partial_t \vc{U} +  \vc{U} \cdot \Grad \vc{U} \Big) \cdot (\vc{U} - \vu) } \ \dt
- \int_0^\tau \intO{ p(\vr, \vt) \Div \vc{U}  } \ \dt
\]
\[
- \int_0^\tau \intO{ \left( \vr \Big( s (\vr, \vt) - s(r, \Theta) \Big) \partial_t \Theta
+ \vr \Big( s(\vr, \vt) - s(r, \Theta) \Big) \vc{U} \cdot \Grad \Theta \right) } \ \dt
\]
\[
+ \int_0^\tau \intO{ \left( \left( 1 - \frac{\vr}{r} \right) \partial_t p(r, \Theta) -
\frac{\vr}{r} \vu \cdot \Grad p(r, \Theta) \right) } \ \dt
\]
for any trio of (smooth) test functions $[r, \Theta, \vc{U}]$ such that
\bFormula{ws4}
r, \ \Theta > 0  \ \mbox{in}\ \Ov{\Omega}, \ \vc{U}|_{\partial \Omega} = 0,
\eF
see \cite{FeiNov10}.

\subsection{Weak-strong uniqueness}

Our goal is to show the following result.

\Cbox{Cgrey}{

\bTheorem{ws1}

In addition to the hypotheses of Theorem \ref{Tg1}, let
\bFormula{bet}
0\le\beta< 1,\;\;\alpha\ge \frac{8}{5(1-\beta)}.
\eF
Suppose that the Navier-Stokes-Fourier system (\ref{i1}--\ref{i4}) admits a strong (classical) solution
$[\tilde \vr, \tilde \vt, \tilde \vu]$ in $[0,T] \times \Omega$. Let $[\vr, \vu, \vt]$ be a weak solution emanating from the same initial data.

Then
\[
\vr \equiv \tilde \vr,\ \vt \equiv \tilde \vt,\ \vu \equiv \tilde \vu \ \mbox{in}\ [0,T] \times \Omega.
\]

\eT

}

The rest of this section is devoted to the proof of Theorem \ref{Tws1}. We point out that similar results proved in \cite{FeiNov10} lean heavily on the assumption that the viscosity coefficients are coercive, in particular, $\mu(\vt) \to \infty$ as $\vt \to \infty$. This is not the case in our setting which makes the proof more delicate. The main idea is to combine the relative energy inequality with the bounds obtained via the thermal energy
balance (\ref{g6}), specifically, (\ref{g32}).

\medskip

\noindent
{\bf Step 1:}

\medskip

We start by introducing the notation:
\[
h = [h ] _{\rm ess} + [h ]_{\rm res},
\]
where
\[
[h]_{\rm ess} = \Phi (\vr , \vt) h , \ \Phi \in \DC ((0, \infty)^2),
\]
\[
0 \leq \Phi \leq 1,\ \Phi = 1
\ \mbox{in an open neighborhood of a compact}\ K \subset (0, \infty)^2
\]
where $K$ is chosen to contain the range of the strong solution $[ \tvr, \tvt ]$, specifically,
\[
[\tilde \vr (t,x) , \tilde \vt (t,x) ] \in K \ \mbox{for all}\ x \in \Ov{\Omega},\ t \in [0,T].
\]

Now, we take $r = \tilde \vr$, $\Theta = \tilde \vt$, $\vc{U} = \tvu$ as test functions in the relative energy inequality
(\ref{ws3}). After a bit tedious manipulation performed in detail in \cite[Section 6]{Feireisl2012} we obtain
\bFormula{ws5}
\left[ \mathcal{E} \left( \vr, \vt, \vu \Big| \tvr, \tvt, \tvu \right) \right]_{t=0}^{t = \tau}
\eF
\[
+ \int_0^\tau \intO{ \left[ \left( \frac{\tvt}{\vt} - 1 \right) \tn{S}(\vt, \Grad \vu) : \Grad \vu + \left( \frac{\vt}{\tvt} - 1 \right) \tilde{\tn{S}} (\tvt, \Grad \tvu) : \Grad
\tvu  \right]  } \ \dt
\]
\[
+  \int_0^\tau \intO{ \left[   \left( \tn{S}(\vt, \Grad \vu) - \tn{S} (\tvt, \Grad \tvu)\right) : \left( \Grad \vu - \Grad \tvu \right) \right]  } \ \dt
\]
\[
+ \int_0^\tau \intO{ \left[ \left( 1 - \frac{\tvt}{\vt} \right) \frac{ \vc{q}(\vt, \Grad \vt) \cdot \Grad \vt }{\vt} + \left( 1 - \frac{\vt}{\tvt} \right) \frac{ {\vc{q}}(\tvt, \Grad \tvt) \cdot \Grad \tvt }{\tvt}
\right] } \ \dt
\]
\[
+ \int_0^\tau \intO{ \left[
\left( \frac{\vc{q}(\vt, \Grad \vt)}{\vt} - \frac{\vc{q} (\tvt, \Grad \tvt) }{ \tilde \vt}   \right) \cdot \left( \Grad \tvt - \Grad \vt \right) \right] } \ \dt
\]
\[
\leq \int_0^\tau \chi \left[ \mathcal{E} \left( \vr, \vt, \vu \Big| \tvr, \tvt, \tvu \right)  + \intO{ \left| \left[p (\vr, \vt)\right]_{\rm res} \right| } \right] \ \dt
\]
\[
\int_0^\tau  \intO{ \Big( \Big[ 1 + \vr + \vr |s(\vr, \vt) | \Big]_{\rm res} \Big) \Big( \Big[ 1 +
\left| \vu - \tvu  \right| \Big]_{\rm res} \Big) + \left| \Big[  \vt - \tvt  \Big]_{\rm res} \right| } \ \dt
\]
where $\chi \in L^1(0,T)$, see \cite[Section 6, formula (71)]{Feireisl2012}.

\medskip

\noindent
{\bf Step 2:}

\medskip

We can write
\[
 \left( \frac{\tvt}{\vt} - 1 \right) \tn{S}(\vt, \Grad \vu)  : \Grad \vu + \left( \frac{\vt}{\tvt} - 1 \right) {\tn{S}}(\tvt, \Grad \tvu) : \Grad
\tvu  +     \left( \tn{S}(\vt, \Grad \vu) - {\tn{S}}(\tvt, \Grad \tvu) \right) : \left( \Grad \vu - \Grad \tvu \right)
\]
\[
= \frac{\tvt}{\vt} \tn{S}(\vt, \Grad \vu) : \Grad \vu + \frac{\vt}{\tvt}{\tn{S}}(\tvt, \Grad \tvu) : \Grad \tvu - \tn{S}(\vt, \Grad \vu) : \Grad \tvu -
{\tn{S}} (\tvt, \Grad \tvu): \Grad \vu
\]
and a similar identity can be derived also for the heat flux.
Thus we have to handle expressions in the form
\bFormula{ws6}
\frac{\tvt}{\vt} \nu |\tn{A}|^2 + \frac{\vt}{\tvt} \tilde \nu | \tilde{\tn{A}} |^2 - (\nu + \tilde \nu) \tn{A} :
\tilde{\tn{A}} = \nu \left| \sqrt{ \frac{ \tvt }{\vt} } \tn{A} \right|^2 + \tilde \nu \left| \sqrt{ \frac{ \vt }{\tvt} } \tilde{\tn{A}} \right|^2 - (\nu + \tilde \nu) \left( \sqrt{ \frac{ \tvt }{\vt} }\tn{A} \right) :
\left( \sqrt{ \frac{ \vt }{\tvt} }\tilde{\tn{A}} \right)
\eF
\[
= \tvt
\frac{\nu}{\vt} \left(   \tn{A}  -  { \frac{ \vt }{\tvt} } \tilde{\tn{A}}  \right)^2
+ (\nu - \tilde \nu) \left( \tn{A} - \frac{\vt}{\tvt} \tilde{ \tn{A} } \right): \tilde {\tn{A} }
\]
\[
\geq \tvt
\frac{\nu}{\vt} \left(   \tn{A}  -  \tilde{\tn{A}}  \right)^2 + 2 \tvt
\frac{\nu}{\vt} \left( 1 - \frac{\vt}{\tvt} \right) \left( \tn{A} - \tilde{\tn{A}} \right): \tilde {\tn{A}}  +
(\nu - \tilde \nu) \left( \tn{A} - \tilde{\tn{A}} \right) : \tilde {\tn{A}} + (\nu - \tilde \nu) \left(1 - \frac{\vt}{\tvt} \right) |\tilde{\tn{A}}|^2,
\]
with the transport coefficient $\nu = \nu(\vt)$, $\tilde \nu = \nu(\tilde \vt)$. If $\nu$ is a continuously differentiable function of $\vt$, we get
\bFormula{ws7}
\intO{ \left[ \frac{\tvt}{\vt} \nu |\tn{A}|^2 + \frac{\vt}{\tvt} \tilde \nu | \tilde{\tn{A}} |^2 - (\nu + \tilde \nu) \tn{A} :
\tilde{\tn{A}} \right]_{\rm ess} }
\eF
\[
\geq c_1 \left\| \left[ \tn{A} - \tilde{ \tn{A} } \right]_{\rm ess} \right\|^2_{L^2(\Omega; R^{3 \times 3})}
- c_2 \intO{  | [ \vt - \tvt ]_{\rm ess} |^2  }
\]
\[
\geq
c_1 \left\| \left[ \tn{A} - \tilde{ \tn{A} } \right]_{\rm ess} \right\|^2_{L^2(\Omega; R^{3 \times 3})}
- c_3 \mathcal{E} \left\{ \vr, \vt , \vu \Big| \tvr, \tvt, \tvu \right\}.
\]

On the other hand, by the standard interpolation argument,
\bFormula{ws8}
\left\| \tn{A} - \tilde{\tn{A}} \right\|_{L^\Lambda(M, R^{3 \times 3})}^2 = \left\| \sqrt{ \frac{\vt}{\nu} } \sqrt{ \frac{\nu}{\vt} }\left( \tn{A} - \tilde{\tn{A}} \right) \right\|_{L^\Lambda(M, R^{3 \times 3})}^2
\eF
\[
\leq
\left\|  { \frac{\vt}{\nu} } \right\|_{L^a(M)} \int_M \frac{\nu}{\vt} \left| \tn{A} - \tilde{ \tn{A} } \right|^2 \ \dx,\
a = \frac{\Lambda}{2 - \Lambda}, \ 1 \leq \Lambda \leq 2
\]
for any measurable set $M \subset \Omega$. Thus for $\Lambda = \frac{8}{5 + \beta}$ (cf. hypothesis (\ref{bet})) we get
\bFormula{ws9}
\left\| \left[ \tn{A} - \tilde{\tn{A}} \right]_{\rm res} \right\|_{L^\Lambda(\Omega; R^{3 \times 3})} \leq \left(
\intO{ \vt^4 } \right)^{(1 + \beta)/4} \intO{ \frac{\mu(\vt)}{\vt} \left| \left[ \tn{A} - \tilde { \tn{A} } \right]_{\rm res} \right|^2 }.
\eF

In accordance with the energy estimates (\ref{g11}), the norm $\| \vt \|_{L^4(\Omega)}$ is uniformly bounded in time for any weak solution, in particular,
relations (\ref{ws7}), (\ref{ws9}) yield
\bFormula{ws10}
\int_0^\tau \intO{ \left[ \left( \frac{\tvt}{\vt} - 1 \right) \tn{S}(\vt, \Grad \vu) : \Grad \vu + \left( \frac{\vt}{\tvt} - 1 \right) \tilde{\tn{S}} (\tvt, \Grad \tvu) : \Grad
\tvu  \right]  } \ \dt
\eF
\[
+  \int_0^\tau \intO{ \left[   \left( \tn{S}(\vt, \Grad \vu) - \tn{S} (\tvt, \Grad \tvu)\right) : \left( \Grad \vu - \Grad \tvu \right) \right]  } \ \dt
\]
\[
\geq c_1 \int_0^\tau \left\| \vu - \tvu \right\|^2_{W^{1,\Lambda}_0 (\Omega;R^3)} \ \dt - c_2 \int_0^\tau \mathcal{E} \left( \vr, \vt , \vu | \tvr, \tvt, \tvu \right) \ \dt , \ \Lambda = \frac{8}{5 + \beta},\ c_1 > 0.
\]

Now, we apply the same arguments to the temperature dependent dissipative terms on the left-hand side of inequality (\ref{ws5}) to obtain:
\bFormula{ws11}
\int_0^\tau \intO{ \left[ \left( 1 - \frac{\tvt}{\vt} \right) \frac{ \vc{q}(\vt, \Grad \vt) \cdot \Grad \vt }{\vt} + \left( 1 - \frac{\vt}{\tvt} \right) \frac{ {\vc{q}}(\tvt, \Grad \tvt) \cdot \Grad \tvt }{\tvt}
\right] } \ \dt
\eF
\[
+ \int_0^\tau \intO{ \left[
\left( \frac{\vc{q}(\vt, \Grad \vt)}{\vt} - \frac{\vc{q} (\tvt, \Grad \tvt) }{ \tilde \vt}   \right) \cdot \left( \Grad \tvt - \Grad \vt \right) \right] } \ \dt
\]
\[
\geq c_1 \int_0^\tau \left[ \left\| \vt - \tvt \right\|^2_{W^{1,2} (\Omega)} + \left\| \Grad (\log(\vt) - \log(\tvt)) \right\|^2_{L^2(\Omega;R^3)} \right]  \dt - c_2 \int_0^\tau \mathcal{E} \left( \vr, \vt , \vu | \tvr, \tvt, \tvu \right) \ \dt,
\]
$c_1 > 0$,
see \cite[Section 6, formula (78)]{Feireisl2012}. Moreover, in view of (\ref{g19}), (\ref{g22}), we have
\bFormula{ws11+}
\intO{ \vr |[\log(\vt)]_{\rm res}| } \leq c \intO{\left([1]_{\rm res} +  \vr |[s(\vr, \vt)]_{\rm res}|+
\vr [e(\vr,\vt)]_{\rm res}  \right)};
\eF
whence, as a consequence of (\ref{ws2}), and the generalized Korn's inequality (cf. \cite[Appendix, Theorem 10.17]{FENO6})
\[
\| \vc{v} \|_{W^{1,2}(\Omega; R^3)} \leq c \left( \| \Grad \vc{v} \|_{L^2(\Omega; R^{3 \times 3})} + \intO{ \vr |\vc{v}| } \right),
\]
the $\Grad \log(\vt)$-norm can be replaced by the $W^{1,2}-$norm.

Going back to (\ref{ws5}) and using (\ref{ws10}), (\ref{ws11}) we may infer that
\bFormula{ws12}
\left[ \mathcal{E} \left( \vr, \vt, \vu \Big| \tvr, \tvt, \tvu \right) \right]_{t=0}^{t = \tau}
\eF
\[
+  c \int_0^\tau \left[ \left\| \vt - \tvt \right\|^2_{W^{1,2} (\Omega)} + \left\| \log(\vt) - \log(\tvt) \right\|^2_{W^{1,2} (\Omega)} + \left\| \vu - \tvu \right\|^2_{W^{1,\Lambda}_0 (\Omega;R^3)} \right]  \ \dt
\]
\[
\leq \int_0^\tau \chi \left[ \mathcal{E} \left( \vr, \vt, \vu \Big| \tvr, \tvt, \tvu \right)  + \intO{ \left| \left[p (\vr, \vt)\right]_{\rm res} \right| } \right] \ \dt
\]
\[
\int_0^\tau  \intO{ \Big( \Big[ 1 + \vr + \vr |s(\vr, \vt) | \Big]_{\rm res} \Big) \Big( \Big[ 1 +
\left| \vu - \tvu  \right| \Big]_{\rm res} \Big) + \left| \Big[  \vt - \tvt  \Big]_{\rm res} \right| } \ \dt
\]
where
\[
\chi \in L^1(0,T),\ \Lambda = \frac{8}{5 + \beta}.
\]

\medskip

\noindent
{\bf Step 3:}

\medskip

In accordance with the energy estimates (\ref{g12}) and the coercivity properties of the relative energy stated in (\ref{ws2}), several terms on the
right-hand side of (\ref{ws12}) are dominated by $\mathcal{E} \left( \vr, \vt, \vu \Big| \tvr, \tvt, \tvu \right)$. Accordingly, (\ref{ws12}) reduces to
\bFormula{ws13}
\left[ \mathcal{E} \left( \vr, \vt, \vu \Big| \tvr, \tvt, \tvu \right) \right]_{t=0}^{t = \tau}
\eF
\[
+  c \int_0^\tau \left[ \left\| \vt - \tvt \right\|^2_{W^{1,2} (\Omega)} + \left\| \log(\vt) - \log(\tvt) \right\|^2_{W^{1,2} (\Omega)} + \left\| \vu - \tvu \right\|^2_{W^{1,\Lambda}_0 (\Omega;R^3)} \right]  \ \dt
\]
\[
\leq \int_0^\tau \chi \mathcal{E} \left( \vr, \vt, \vu \Big| \tvr, \tvt, \tvu \right)  \dt
+
\int_0^\tau  \intO{ \Big[ \vr |s(\vr, \vt) | \Big]_{\rm res}  \Big[
\left| \vu - \tvu  \right| \Big]_{\rm res}   } \ \dt.
\]

Consequently, our ultimate goal is to ``absorb'' the last integral on the right-hand side of (\ref{ws13}) via a Gronwall type argument. In accordance with
(\ref{g18}), (\ref{g19}) we have to estimates the integral
\[
\int_0^\tau \intO{ \left[ \vr^{\gamma/3} + \vr \vt^\omega + \vr |\log \vt| + \vt^3 \right]_{\rm res} \left[ \vu - \tvu \right]_{\rm res} } \ \dt.
\]

First, we observe that
\[
\intO{ \left[ \vr^{\gamma/3} + \vr \vt^\omega  \right]_{\rm res} \left[ \vu - \tvu \right]_{\rm res} }
\]
\[
\leq \frac{1}{2} \intO{ \vr |\vu - \tvu |^2 } + \frac{1}{2} \intO{ \left[ \vr \vt^{ 2\omega} + \vr^{1 + 2 \gamma /3 } \right]_{\rm res} }.
\]
As $0 \leq \omega \leq 1/2$ and $\gamma \geq 3$, the second integral is dominated by $\mathcal{E} \left( \vr, \vt, \vu \Big| \tvr, \tvt, \tvu \right)$.

Next, we estimate
\[
\intO{ \left[ \vr |\log (\vt)| \right]_{\rm res} |\vu - \tvu | } \leq \left\| \sqrt{\vr} \right\|_{L^{3}(\Omega)}
\left\| \left[ \log (\vt) \right]_{\rm res} \right\|_{L^6(\Omega)} \left\| \sqrt{\vr} (\vu - \tvu) \right\|_{L^2(\Omega;R^3)}
\]
\[
\leq c\left( \delta \left\| \left[ \log(\vt) \right]_{\rm res} \right\|_{W^{1,2}(\Omega)}^2 + c(\delta) \intO{ \vr |\vu - \tvu |^2 }\right),
\]
where the last integral is dominated by the terms on the left-hand side of (\ref{ws13}) for $\delta > 0$ small enough.

Finally, by H\" older's inequality,
\[
\intO{ \left[ \vt^3 \right]_{\rm res} | \vu - \tvu | }
 \leq \| \vu - \tvu \|_{L^a(\Omega; R^3)} \Big(\intO{[\vt^4]_{\rm res}}\Big)^{1/2} \,\Big(\intO{\vt^b}\Big)^{1/b}
\]
\[
\leq \delta \| \vu - \tvu \|_{L^a(\Omega; R^3)}^2 + c(\delta)\|\vt\|_{L^b(\Omega)}^2 \mathcal{E} \left( \vr, \vt, \vu \Big| \tvr, \tvt, \tvu \right)
\]
where $a=\frac{24}{7+5\beta}$ is determined from the Sobolev embedding $W^{1,\Lambda}\subset L^a$, $\delta>0$ is arbitrary, $c(\delta)>0$, $\frac1a+\frac12+\frac1b=1$ and $\|\vt\|_{L^b(\Omega)}^2\in L^1(0,T)$ by virtue
of (\ref{g25}) provided
$$
\alpha\ge \frac 8{5(1-\beta)},\; 0\le \beta< 1.
$$
We have proved Theorem \ref{Tws1}.

\section{Conditional regularity}\label{cr}

Similar to \cite{FeNoSun1}, the weak-strong uniqueness principle established in Theorem \ref{Tws1}
can be used to deduce a \emph{regularity criterion} in the class of weak solutions:

\Cbox{Cgrey}{

\bTheorem{cr1}
In addition to the structural restrictions imposed in Theorem \ref{Tws1} on the transport coefficients $\mu$, $\eta$, and $\kappa$, assume
that there exists a constant $k>0$ such that
\bFormula{cr2}
\mu(\vt) = k \eta(\vt).
\eF
Let $[\vr,\vu,\vt]$ be a weak solution to the Navier-Stokes-Fourier system (\ref{i1})-(\ref{i3}) in the time interval $(0,T)$, satisfying
the boundary conditions (\ref{i4}), the initial conditions (\ref{g1}), where
\bFormula{cr3}
\vr_0, \, \vt_0 \in W^{3,2}(\Omega), \, \vu \in W^{3,2}(\Omega; R^3),
\vr_0 \ge \underline{\vt} >0, \, \vt_0 \ge \underline{\vt} >0,
\eF
with suitable compatibility conditions. If
\bFormula{cr4}
{\rm ess}\sup_{t\in (0,T)}\|\Grad\vu(t,\cdot)\|_{L^\infty(\Omega; R^{3\times 3})} < \infty,
\eF
then $[\vr,\vu,\vt]$ is in fact a classical solution.

If, in addition, $\beta=0$, meaning $\mu$, $\eta$ are positive constants independent of $\vt$, then the regularity condition (\ref{cr4}) can be relaxed to
\bFormula{cr5}
\int_0^T \| \Grad\vu(t,\cdot) \|^2_{L^\infty(\Omega; R^{3\times 3})} \ \dt < \infty.
\eF

\eT

}

In view of the weak-strong uniqueness principle, we only need to establish the following \emph{blow up criterion} for a strong solution $[\tvr,\tvu,\tvt]$ with the initial data $[\vr_0,\vu_0,\vt_0]$. We note that the existence of (local) strong solution is guaranteed by
the result of Valli  \cite{Vall1}, \cite{Vall2}.

\Cbox{Cgrey}{

\bProposition{cr7}
Given $[\vr_0,\vu_0,\vt_0]$ satisfying (\ref{cr3}), suppose $[\tvr,\tvu,\tvt]$ is a strong solution to the initial boundary value problem for the Navier-Stokes-Fourier system (\ref{i1})-(\ref{i3}) in the time interval $(0,\tilde{T})$, such that
\[
\tvr,\tvt \in C([0,\tilde{T}), W^{3,2}(\Omega)), \, \tvu \in C([0,\tilde{T}), W^{3,2}(\Omega; R^3)).
\]
If (\ref{cr4}) (or (\ref{cr5}) if $\mu,\eta$ are constants) holds with $\vu$ replaced by $\tvu$, then $[\tvr,\tvu,\tvt]$ can be continuously extended beyond the time $\tilde{T}$.
\eP

}

To prove Proposition \ref{Pcr7} we can follow step by step the proof in \cite{FeNoSun1}, Section 4, keeping in mind the following observations.

1. In \cite{FeNoSun1} we assume $\eta=0$. However this does not hold anymore in our present setting. Condition (\ref{cr2}) ensures that the Lam\'e operator
\[
\mu(\vt) \Delta_x \vu + \left( \eta(\vt) + \frac{1}{3}\mu(\vt) \right)\Grad \Div\vu
\]
to be an elliptic operator with constant coefficients after divided by $\eta(\vt)$. See also remarks in the final section in \cite{FeNoSun1}.

2. Note that assumption (\ref{cr4}) is only used to control the divergence part of the term $\vc{h}_1$ in order to yield estimates (4.15) in \cite{FeNoSun1}, due to the dependence of $\mu$ on $\vt$. If $\mu$(and $\eta$) is a constant, then this divergence term disappears and one can replace (\ref{cr4}) with (\ref{cr5}).

\def\cprime{$'$} \def\ocirc#1{\ifmmode\setbox0=\hbox{$#1$}\dimen0=\ht0
  \advance\dimen0 by1pt\rlap{\hbox to\wd0{\hss\raise\dimen0
  \hbox{\hskip.2em$\scriptscriptstyle\circ$}\hss}}#1\else {\accent"17 #1}\fi}

\end{document}